\documentclass[11pt]{article}
\pdfoutput=1
\usepackage{graphicx}
\usepackage{amssymb}
\usepackage{epstopdf}
\usepackage{amsmath}
\renewcommand{\author}[1]{\large\rm #1\\ \bigskip}
\newcommand{\address}[1]{{\normalsize\it #1\\}\bigskip}
\renewcommand{\title}[1]{\bigskip\bigskip\Large\bf #1\bigskip\bigskip\\}

\begin{document}
	\begin{center}
		\title{The design of efficient algorithms for enumeration}
		\author{Andrew R. Conway\footnote[1]{email:  {\tt andrewenumeration@greatcactus.org}  }  }
		
		\address{ ARC Centre of Excellence for Mathematics and Statistics of Complex Systems,\\
			Department of Mathematics and Statistics,\\
			The University of Melbourne, Victoria 3010, Australia}
	\end{center}
	\begin{abstract}
	Many algorithms have been developed for enumerating various combinatorial objects in time exponentially less than the number of objects. Two common classes of algorithms are dynamic programming and the transfer matrix method. 
	This paper covers the design and implementation of such algorithms.
	
	A host of general techniques for improving efficiency are described. Three quite different example problems are used for examples: 1324 pattern avoiding permutations, three-dimensional polycubes, and two-dimensional directed animals.
	
	For those new to the field, this paper is designed to be an introduction to many of the tricks for producing efficient enumeration algorithms. For those more
	experienced, it will hopefully help them understand the interrelationship and implications of a variety of techniques, many or most of which will be familiar. 
	The author certainly found his understanding improved as a result of writing this paper.
	\end{abstract}
	
\section{Introduction}
	
There are many problems in combinatorics where one wants to know the number of some type of object for a given size $n$ (e.g the number of polygons of perimeter $n$ on a given lattice\cite{enting1980sap}),
but no formula is known. In this case computer enumeration is helpful to get the first few terms of the series, which can then be analysed to estimate their
asymptotic behaviour. For many combinatorial problems the series terms
grow rapidly - frequently exponentially - and it rapidly becomes computationally intractable to examine each object individually. An algorithm that groups multiple objects
together can be more efficient, and can enable enumeration of more terms, which is usually highly desirable. Many such algorithms have been designed for many different
problems. Many of them can be described as dynamic programming algorithms or as transfer matrix algorithms. This paper gives ways ot thinking about and optimizing such algorithms.

The paper does not intend to be exhaustive - that would be impractical - but rather to describe a coherent group of techniques in a common language.
	
\subsection{Enumeration problems as a set of equations}

Ignoring efficiency for the moment, consider how one could count a combinatorial object of size $n$. Often this involves starting with some base state $S_b$ (which typically includes $n$)
and adding something to it - in the case of a geometrical object, this could be the presence or absence of a bond or a site. 
This then produces a set of other states $S_i$ which consist of the base state plus the information about
the new site. For each of these states $S_i$, one adds the next possible item to get a new set of (more complex) states. Continue doing this recursively until one gets to an end state, at which point
one increments a counter and then continues with the next substate (or {\it child} state) in the recursive algorithm. At the end the answer is what is left in the counter. 

\begin{figure}
	\centering
	\includegraphics[width=\linewidth]{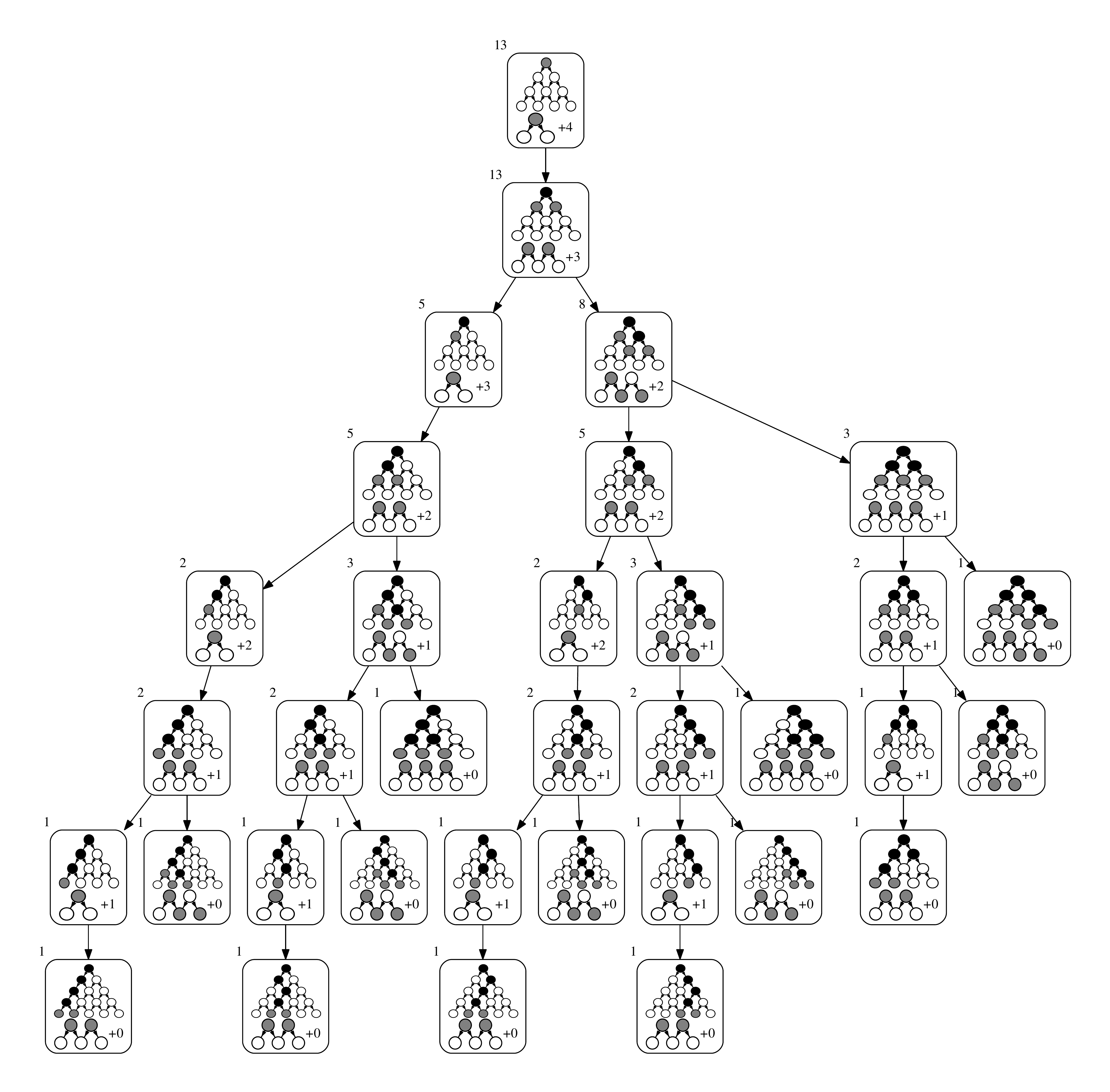}
	\caption{Generation of all directed animals of 4 steps. Each box contains, at top, the animal so far. Black sites are occupied; grey sites are potentially occupied, and white
		sites are unoccupied or not reached yet. At each step, the rightmost grey site in the uppermost row containing a grey site is processed. It is determined to
		be either unoccupied (left child) or occupied (right child). If it becomes occupied, then the sites accessible from it are set to grey. If the site is the last
		grey, then it may not be unoccupied as such a child would never be able to grow to the desired size, so there is only one child. If the animal is the desired
		size, it is finished and has no children. Below the picture of the animal so far is derived information on the border; the grey sites (ignoring surrounding, unreachable sites),
		and the number of sites left to be occupied to reach the desired size. Above and to the left of each box is the number of directed animals descending from that
		state; it is one for a leaf box and otherwise equal to the sum of its direct children. The root node's number, 13, is the number of four site directed animals.}
	\label{fig:FullDirAnTree}
\end{figure}

The set of states and connections to their child states defines a directed acyclic graph, called the {\it call graph} as it represents the calls in the recursive
evaluation of the set of equations. Example call graphs for different formulations of the same problem are given in figures \ref{fig:FullDirAnTree}, \ref{fig:PartialDirAnTree} and
\ref{fig:PartialDirAnTree2}.

This can be written out more formally as a set of equations. Let $f(S)$ be the number of objects with state $S$. So what one wants to calculate is $f(S_b)$. For each state $S$ define
a set of next-possible-states (or {\it child} states) 
$N(S)$ as the states one can get to by doing the next thing to $S$. This may be empty for an end state. Let $k(S)$ be a constant (typically 0, or 1 if
$S$ is an end state). Then one can say
$$f(S)=k(S)+\sum_{s\in N(S)}f(s).$$

We require that there be no loops in the directed graph implied by $N(S)$ in order to make the above equation, in principle, straightforward to evaluate without having to solve simultaneous equations. 
That is, there is no set of $k>1$ states ${s_i}$, such that $s_{i+1}\in N(s_i)$ and $s_1=s_k$. This is normally obvious, and one can often define a {\it hierarchy} function from states to integers (or other
well ordered values), $h(s)$, with
the property $\forall s,\forall t\in N(s), h(t)<h(s)$. Clearly if there exists a hierarchy function, there cannot be any loops. A hierarchy function is called {\it ideal} if $\forall s,\forall t\in N(s), h(t)=h(s)-1$.

An example is given in figure \ref{fig:FullDirAnTree} which builds up all directed animals on a square lattice of size 4. Each rectangle is a state, with arrows pointing to its child states. 
The hierarchy function of a state in this case is a pair of integers; the number of sites left to add, and the number of grey sites. Ordering on the tuples is primarily by the number of sites to
add, and secondarily by the number of grey sites.

Frequently an ideal hierarchy function is obvious --- the number of elements to add to the
partially constructed object, or the number of sites left to consider when passing over a finite lattice.

Sometimes there are weights associated with the substates - this is particularly common when the value of $f$ is a partial generating function (typically of integers) rather than an integer.
Let $w(s,S)$ be the (easy to compute) weight associated with state $s$ coming from state $S$, giving us a more general formula
$$f(S)=k(S)+\sum_{s\in N(S)}w(s,S)f(s)$$

Sometimes the series term number, $n$, is not explicitly in the state but is implicitly stored and is represented in the structure of $f$ (e.g. the size of the lattice when
enumerating on a finite lattice) and/or in the size of the partial generating function that is the value of $f$.

We call a definition of $N$ {\it clean} if there is no path from the base state to a state $s$ where $f(s)=0$. If $N$ is not clean then extra, useless, work will be done. Unfortunately
it can be difficult to be clean.

It turns out that a large class of combinatorial enumeration problems fall into this formalism, although the notation may be somewhat different. Sometimes the formalism
is slightly different, requiring modifications of the techniques below (e.g. section \ref{sec:Factorizations})

As an example, consider enumerating the set of well-formed bracket expressions with $n$ open and $n$ closing brackets\footnote{A well formed bracket(or parenthesis) expression has the same number of 
open and closing brackets, and there are never more closing brackets than open brackets for any prefix of the expression. So ()() is OK, as is (()), but not ())(. This is a well known
(and solved) problem whose series is the Catalan numbers.}.
Define a state as a pair of integers, $S=(o,c)$, meaning the state where one has $o$ open and $c$ closing brackets left to add. So the stating state would be $S_b=(n,n)$, the ending state
would be $(0,0)$, and adding an open or closing bracket would decrease $o$ or $c$ respectively. Overall one has
$$f(o,c)=
 \left\{
   \begin{array}{ll}
     1 & \mbox{if } o=0 \mbox{ and } c=0 \\
     0 & \mbox{if } o>c \\ 
     f(o-1,c) & \mbox{if } o=c>0 \\
     f(o-1,c)+f(o,c-1) & \mbox{otherwise}
\end{array}
\right.
$$
Note that this definition is redundant. The second line could be omitted, as such states will never be reached as the definition is clean. 
Alternatively, the third line could be omitted, but then the system would not be clean and the second line would be used.
Also note that there are many other ways of enumerating well formed bracket expressions, with different state definitions. 

The straightforward recursive implementation of this function typically takes time proportional to $f(S_b)$ as the only terminal state is $1$.
In the case of the well-formed bracket expressions (and many others of interest), this grows exponentially with $n$.

\subsection{Dynamic programming}

Recursive functions of this kind often have a simple optimization - whenever one has calculated $f(S)$ for some $S$, one adds $S\rightarrow f(S)$ in some table. 
Whenever one try to compute $f(s)$, one checks whether $s$ is present in the table already, and if so use the stored value rather than recomputing it. This is generally called
{\it dynamic programming} (DP),  {\it memoization}, or {\it caching} 

DP takes up memory, but can vastly reduce the running time. Of course it relies on the same state being reachable multiple ways - otherwise the table
will never be used.

As an easy to analyze case study, consider the well formed bracket expressions from above, and suppose one is trying to compute $f(10,10)$.
The first step is forced, $f(10,10)=f(9,10)$. The next step has two legal choices, $f(9,10)=f(8,10)+f(9,9)$. Now comes the interesting bit.
$f(8,10)=f(7,10)+f(8,9)$ and $f(9,9)=f(8,9)$. Both use $f(8,9)$. Using dynamic programming, the second time the program attempts to evaluate
$f(8,9)$ it would just use the value previously computed. 
The total number of legal signatures needed to compute $f(n,n)$ is the number of pairs of integers $(o,c)$ where $n\geq c \geq o \geq 0$. This
is clearly $O(n^2)$ which means the total number of states visited (and thus time and memory use) of the algorithm is $O(n^2)$. This compares
with the naive recursive function which is $O(4^n n^{-1.5})$ which is clearly much worse.

In this case dynamic programming reduced the complexity of the algorithm from exponential to polynomial, and this sometimes happens in
practice. When it does, the problem becomes computationally easy, a long series is generated, and people don't bother enumerating
it again! Frequently indeed it means that it is plausible to get a closed form solution. 
In other problems DP only reduces the exponential growth rate, which can significantly extend the number of terms one can get, although
it rapidly becomes computationally intractable again. People tend to spend most of their time on such problems, as they are the unsolved ones.

\subsection{Transfer matrix}

Another approach is to keep a set of states and associated values (or {\it multiples}), starting off with $S_b$, multiple $1$. Then, in each iteration of the
algorithm one takes each state $S$ and associated value $v$, and applies the function $f$ to it, adding each of the states $s \in N(S)$ with multiple
$m$ (or $m \cdot w(s,S)$ if weights are present) to the new set of states and associated multiples. 
Also one keeps a counter and whenever applying $f$ to a state $S$ produces a constant, one adds the
associated multiple times the constant to the counter. Then one replaces the old set of states and multiples with the new set of states and multiples and 
iterates until the set of states becomes empty.

So far nothing has been gained; one is effectively doing a breadth first search of the recursive algorithm rather than the obvious depth first
search. But crucially, when one ends up putting the same state into the new set of states multiple times, one can merge them together and
sum the associated multiples. Now, in a very similar manner to DP, one has bunched together a variety of ways to get to the same state, preventing
redundant computation. This technique is often called the {\it transfer matrix technique} (TM).

Again going through the well formed brackets example, we have
\begin{equation}
\begin{aligned}
f(10,10) &= f(9,10) \\
         &= f(8,10)+f(9,9) \\
         &= f(7,10)+2f(8,9) \\
         &= f(6,10)+3f(7,9)+2f(8,8) \\
         &= \dots
\end{aligned}
\end{equation}

Typically a TM algorithm reduces the total amount of computation to the number of valid signatures reachable from the initial state, the same as DP\footnote{If one
can reach the same state via two different routes {\it of different lengths} through $f$, then one will end up with the same state in multiple iterations
and the amount of computation will be higher - worse than DP. This can usually be avoided. In particular, if one has a hierarchy function $h$ defined on states,
then one can, at each pass, process only the states with the highest hierarchy value.}.
The memory use is somewhat less than DP, as one is effectively ordering computation so that one has some states that one can discard information
about once one knows that one can never revisit them. The memory use is the maximum number of states one gets in a breadth first search of
the structure of $f$\footnote{Actually roughly twice this value, as one needs to have two sets of states in memory at once; the set being 
processed and the set being produced}. This is bounded above by the number of states from DP, and bounded below by said number of states divided by
the number of iterations of $f$ before reaching a final state, which is typically $n$ or at least polynomial in $n$. This means that the
memory use is better than DP, but still growing with the same exponent if it is exponential.	
	
\subsection{Comparison and interconversion of DP and TM}

The computational complexity of DP and TM are generally the same, but TM generally has somewhat lower memory requirements. So why would anyone ever use DP?

One reason is algorithm design complexity. TM is somewhat more complex to think about, implement, and test and debug than DP, while useful tools like Maple
sometimes have direct support for DP (Maple calls it memoization). The resources used to implement TM could possibly be better spent making other changes
to the algorithm that give a better return on the implementer's time. 

Another reason is lack of awareness of the trick of doing TM on DP problems.
The main purpose of this paper is to raise awareness of such tricks, and introduce a way of thinking about the
algorithms such that the interconversion is conceptually straightforward.

A third reason is that DP can be used to solve non-linear equations for $f$ - see section \ref{sec:Factorizations}. This can be the real deciding factor.

A TM algorithm can almost always be converted to a DP algorithm, possibly requiring adding some information to the state - 
undoing the effect of section \ref{sec:TMItNoSig}. But there is rarely a point in doing so.

\subsection{Finite Lattice Method}

For enumeration of geometrical objects on a lattice, it is common to split up the problem into enumeration on each possible finite
lattice. Then the TM method is used, where each step consists of dealing with one extra element of the lattice (usually site). The ideal
hierarchy function is the number of elements left to be added. The desired series is then produced by combining the series from the
finite lattices. This is commonly referred to as the Finite Lattice Method. 

\begin{figure}
	\centering
	\includegraphics[scale=1.0]{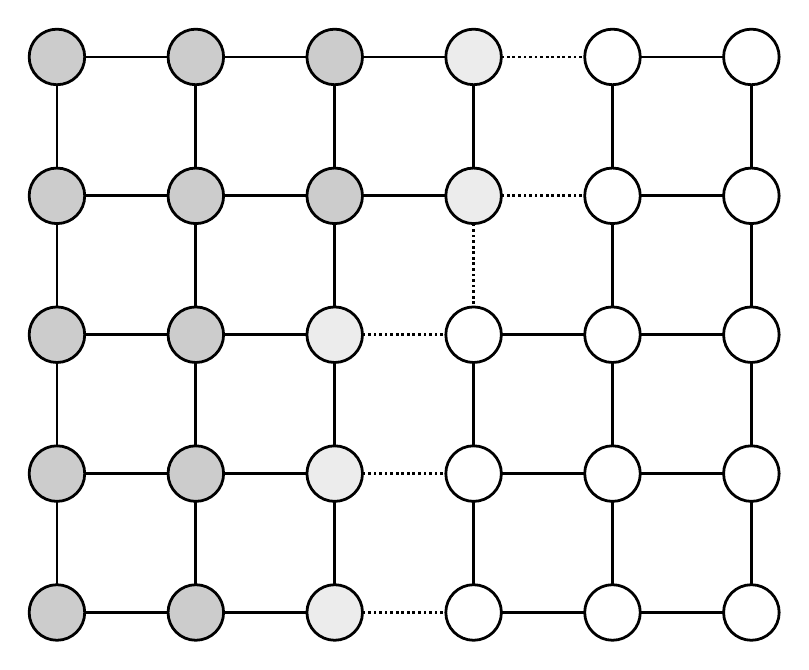}
	\caption{The boundary of a finite lattice. The grey sites are processed, the white sites unprocessed. The light grey sites
		are processed sites on the boundary, the dark grey ones are sites whose status is now irrelevant other than in so far as
		they contribute to connectivity information on the boundary. The boundary may be sites, in which case it is usually the
		grey sites, or bonds, in which case it is usually the dotted bonds. Note the kink as sites are processed one at a time; the
		next site to be processed is the white site adjacent to two grey sites.   }
	\label{fig:2Dboundary}
\end{figure}

Typically the sites are added in a systematic way that minimises the size of the boundary between processed and unprocessed sites,
and the state then describes that boundary and connectivity information. On a rectangular lattice\footnote{A rectangular subset of the square lattice $\mathbb{Z}^2$}
the boundary would be a cross sectional line through the lattice, possibly with a kink in it if the number of sites processed is not an even
multiple of the lattice width. This is shown in figure \ref{fig:2Dboundary}.

Doing it on multiple lattices has several advantages
\begin{itemize}
	\item It makes it easy to deal with uniqueness under translation invariance
	\item One can commonly require that the object on the finite lattice touches all edges of the lattice. While keeping track of the information
	      required to track this increases the total number of legal states, in practice it often decreases the number of states that will be
	      visited in a clean implementation restricted to a particular size.
    \item One can use lattice symmetries. Enumerating on the $n$ by $m$ lattice is usually the same as the $m$ by $n$ lattice.
          The orientation is chosen to minimize the length of the boundary, and thus the number of states.
\end{itemize}

The finite lattice method has been used to great effect on two dimensional lattices, e.g. \cite{conway1995polyonimo,jensen2001polyonimo,jensen2003polyonimo,conway1993saw,conway1996saw,jensen2013saw,enting1980sap,jensen1999sap,jensen2003parallelsap,clisby2012sap}. 
When enumerating size $n$ objects, one generally has a maximum boundary length of roughly $n/2$. The number of states usually grows exponentially
in the boundary length, although pruning can mean that the peak number of states actually encountered is highest at a smaller portion of $n$.

On $d>2$ dimensional lattices, the boundary is now a $d-1$ dimensional surface,
typically of size $(n/d)^{d-1}$. The number of states can grow exponentially in this cross section, which can, in a unclean system,
easily produce a number of states far greater than the number of objects being enumerated, making it a dreadfully inefficient algorithm.
Making a clean system can be difficult, and even if it is clean the performance advantage over direct enumeration is much smaller
than in two dimensions.

A technique occasionally used in the finite lattice method (e.g. \cite{conway1993saw}) is to enumerate subportions of the objects being enumerated and
construct the objects being enumerated from the subportions. This enables smaller lattices to be used, producing smaller
boundary sizes and vast reductions in the number of states.

This has been used for a long time \cite{deNeef1977finite}.

\subsection{Examples}

There are three significantly different examples that will be described in detail. They are used throughout the rest of this paper to provide specific examples of
some techniques and trade-offs.
	
\subsubsection{2D directed animals}

A directed animal on some directed lattice is a set of sites such that each site apart from a single root site is directly downstream of another site.
We will deal here with directed animals on a square lattice $\mathbb{Z}^2$ with the direction constraint being in the positive direction of either the $x$ or $y$ axes (north or east)
Often the lattice will be shown rotated clockwise by 135 degree, so that the directed constraint becomes down and left or down and right.
This is a well solved problem \cite{dhar1982dirAns,gouyou1988dirAns,betrema1993dirAns}; it is included as it is a simple example. Similar
techniques can be used in many other two dimensional problems and on other lattices. It is solved 
for the square and triangular lattice\footnote{The generating
	function $f(x)$ for the square lattice satisfies $(3x-1)f^2+(3x-1)f+x=0$, the triangular similar.}, but not for the hexagonal
or most other lattices. Hereafter we will assume the square lattice.

There is one size one directed animal, the root. There are two size two directed animals, the root and one of the two downstream sites. There are
five size three directed animals: four coming from the root, one of the two downstream sites, and also one of the two further downstream sites; the fifth
is the root and both of the two downstream sites. There are thirteen size four directed animals; their construction is show in figure \ref{fig:FullDirAnTree}. 

\begin{figure}
	\centering
	\includegraphics[scale=0.85]{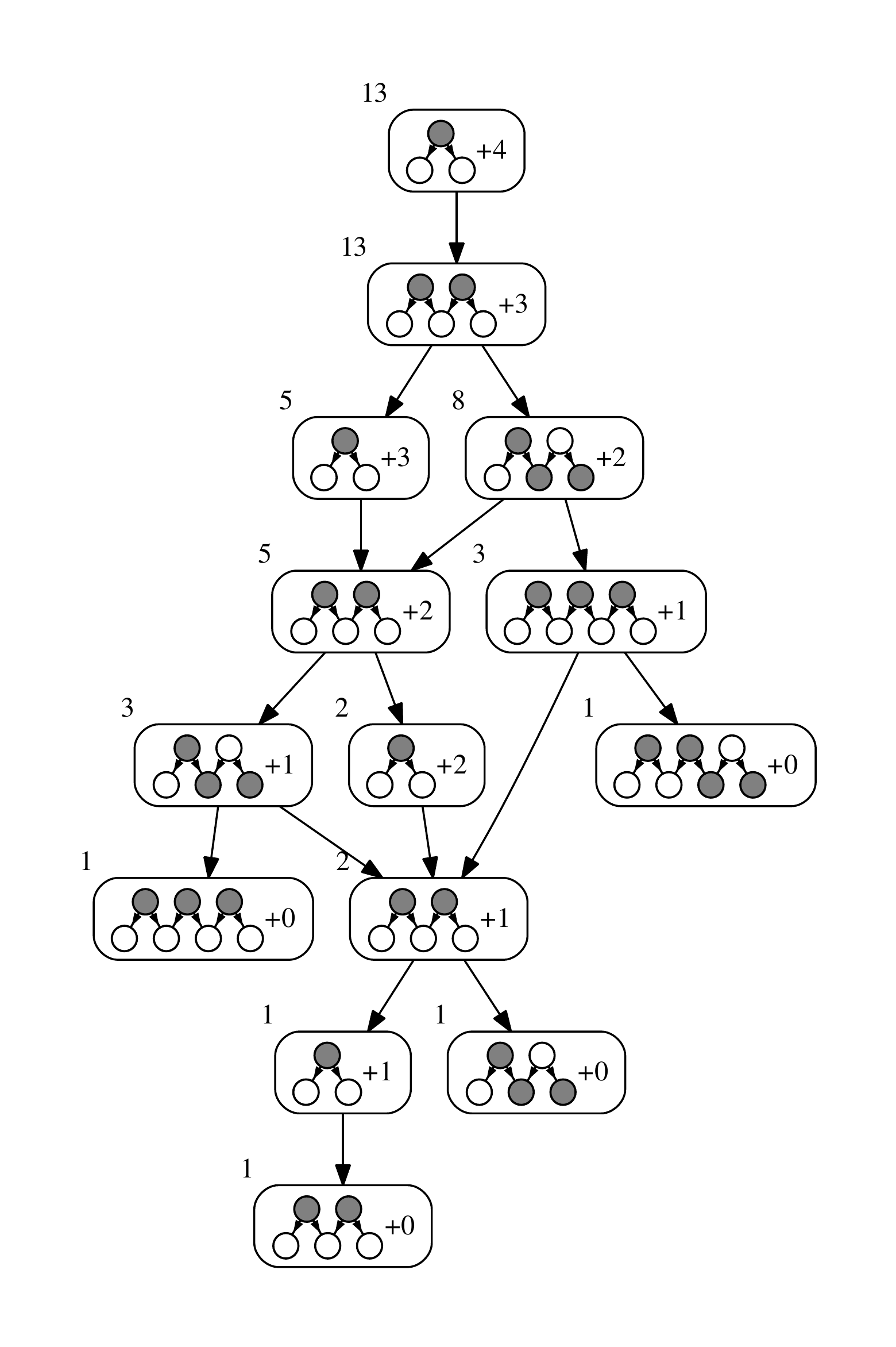}
	\caption{A simplification of figure \ref{fig:FullDirAnTree} just showing the border information, which is the only thing that can affect children. This
		allows coalescing of now identical states together to produce a smaller tree, which will be able to take advantage of dynamic programming. Note there is
		an arrow skipping the fifth row. }
	\label{fig:PartialDirAnTree}
\end{figure}

\begin{figure}
	\centering
	\includegraphics[width=\linewidth]{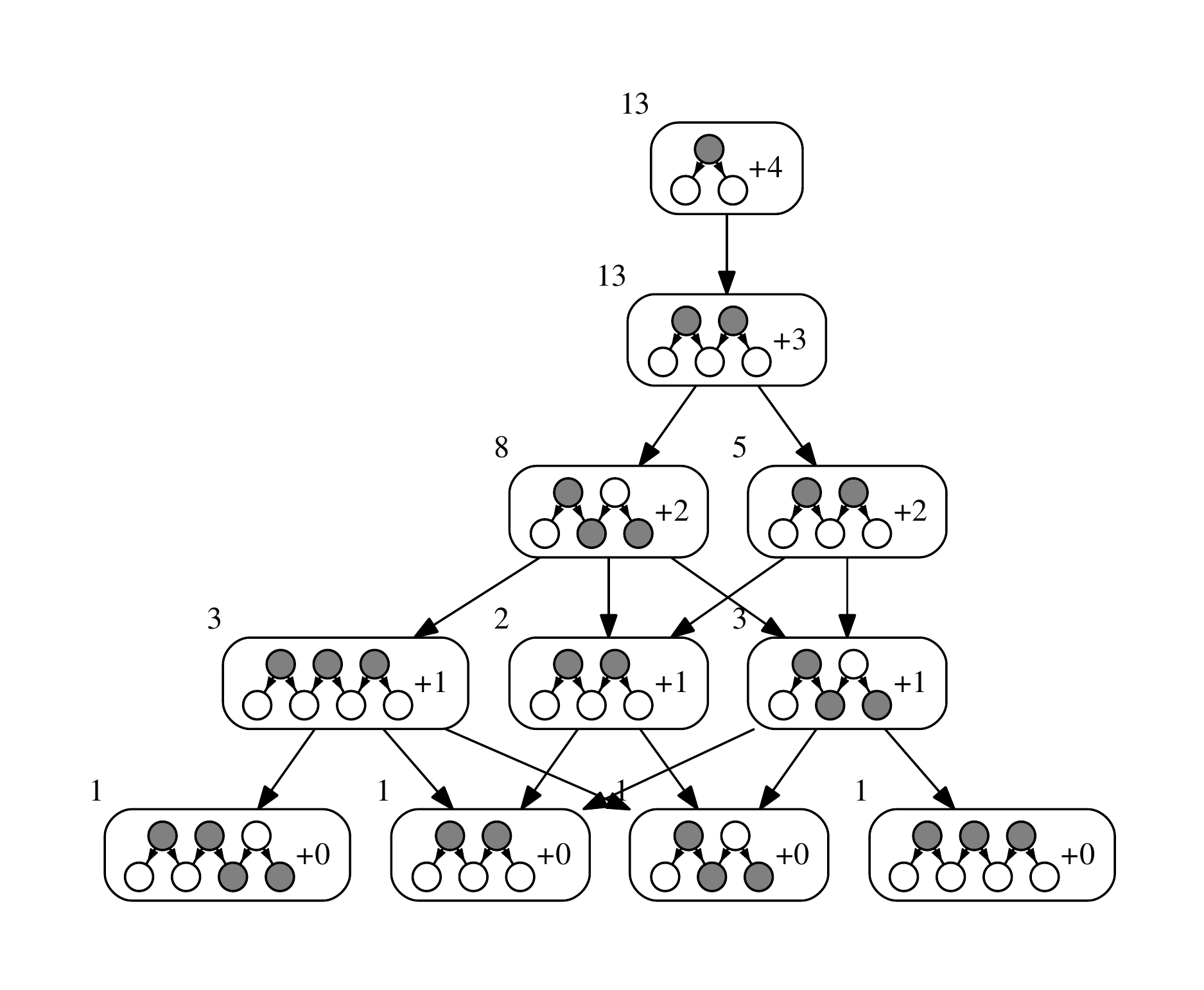}
	\caption{An alternative algorithm for enumeration, like figure \ref{fig:PartialDirAnTree}, except instead of processing one grey site at each step, the algorithm goes through grey sites
		until one is determined to be occupied. This is a more complex algorithm, but is better suited to the transfer matrix method as there are no row skipping arrows - see section \ref{sec:TMItNoSig}. }
	\label{fig:PartialDirAnTree2}
\end{figure}

These can be solved in a straightforward manner using both the dynamic programming and transfer matrix techniques, although the dynamic programming method is more obvious.
One starts from the root, and add sites along a line perpendicular to the preferred direction (north east). Connectivity is then not an issue, as everything reached on this
line must have come from the root, meaning the presence or absence of the sites on this line, along with the number of sites used, is sufficient information. The line
will end up having a kink in it at the site being added, and will be partly the line being processed and partly the subsequent line. In practice one can do even better by
using the (kinked) line one beyond the line to which sites have just been added, and record which sites are reachable rather than which sites have been reached. This is because there can be
multiple sets of reached sites that can produce the same set of reachable sites. The call graph of such a formulation is given in figure \ref{fig:PartialDirAnTree}, with
corresponding equations in table \ref{tab:PartialDirAnTree}. An alternate formulation is given in figure \ref{fig:PartialDirAnTree2} and the transfer matrix algorithm
applied to it is presented in table \ref{tab:TMDirAn}.

\begin{table}
	\begin{tabular}{rcl}
		\hline
		$f(?,?,0)$ & = & $1$, otherwise \\
		$f(1,0,n)$ & = & $f(11,0,n-1)$ \\
		$f(11,0,n)$ & = & $f(1,0,n)+f(10,11,n-1)$ \\
		$f(10,11,n)$ & = & $f(11,0,n)+f(111,0,n-1)$ \\
		$f(111,0,n)$ & = & $f(11,0,n)+f(110,11,n)$ \\
		\hline
	\end{tabular}
	\caption{The equations corresponding to the graph in figure \ref{fig:PartialDirAnTree}. The first two arguments to the function are binary numbers representing
		the pattern of grey sites, with a 1 meaning a grey site. Leading zeros are not shown. The third argument is the number of sites remaining to be added. The three
		arguments together define the state.}
	\label{tab:PartialDirAnTree}
\end{table}

\begin{table}
	\begin{tabular}{rcl}
		\hline
		$f(1,0,4)$ & = & $f(11,0,3)$ \\
		           & = & $f(10,11,2)+f(11,0,2)$ \\
		           & = & $f(111,0,1)+2f(11,0,1)+2f(10,11,1)$ \\
		           & = & $f(110,11,0)+5f(11,0,0)+5f(10,11,0)+2f(111,0,0)$ \\
		           & = & $1+5+5+2$ \\
		           & = & $13$ \\
		\hline
	\end{tabular}
	\caption{Transfer matrix style evaluation of the graph in figure \ref{fig:PartialDirAnTree2}, with function arguments as described in table \ref{tab:PartialDirAnTree}.}
	\label{tab:TMDirAn}
\end{table}

\subsubsection{1324 Pattern avoiding permutations}

A permutation $P$ of the integers $1 \dots n$ is said to avoid the pattern $p$ where $p$ is itself a permutation of $1 \dots l$ if there does not exist an $l$ length subsequence (consecutive or not) of $P$ 
that has the same relative order of elements as in $p$. For instance the permutation $15342$ contains the pattern $123$ through the subsequence $134$.
The enumeration problem is to compute the number of such permutations for a given $n$. 

Surprisingly, all pattern avoiding permutations grow roughly exponentially, rather than the factorial growth of the permutations.\cite{marcus2004excluded}

The series are well understood for all 2, 3, and 4 length patterns except for $1324$ and its complement $4231$. This has been enumerated using a dynamic
programming algorithm. Conceptually the permutations are built up one element at a time, and the state consists of the remaining integers and constraints upon their order. 
This was done in \cite{johansson20141324PAPs} enumerating to 31 terms and improved on in \cite{conway20151324PAPs} with a more efficient state definition extending the series to 36 terms.
The latter algorithm will be used in this paper.

The state definition consists of a series of numbers separated by brackets like $4(5(6)3)5$. Each number represents a set of consecutive integers
yet to be chosen, and the value is the length of this set. The sum of all the numbers is the number of elements left to add. This
sum is an ideal hierarchy function. One may not use an element inside a bracketed expression
until all numbers after the bracketed expression have been dealt with. That is, in the example, there are 4 low numbers and 5 high numbers that are legal as the
next number to be chosen for the permutation being built up. That is, $N$ of the example will have 9 elements. 
There may also be a comma at the start of the state indicating that at least one number lower 
than all remaining numbers has been already taken. See \cite{conway20151324PAPs} for a more precise definition and rules for computing $N$. The start state is
the simple integer $n$, and the end state is the empty state. The system is clean.

\subsubsection{3D polycubes}

\begin{figure}
	\centering
	\includegraphics[scale=1.0]{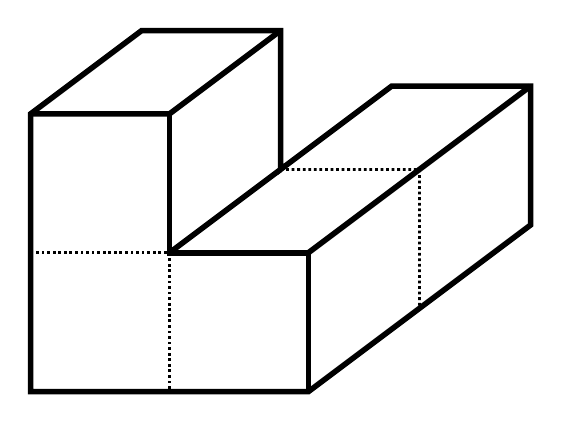}
	\caption{A size (volume) 4 polycube}
	\label{fig:polycube}
\end{figure}

A polycube is a connected set of unit length cubes in $\mathbb{Z}^3$. An example is in figure \ref{fig:polycube}. 
This is the 3D equivalent of a polyonimo. The objective is determining the number of polycubes consisting of $n$ cubes.
This is an open problem in two or higher dimensions, with a roughly exponential growth rate.

There is one size 1 polycube - the single cube. There are three size two polycubes - two adjacent blocks, pointing along one of the three axes.

Polycubes are difficult to count efficiently. There has been extensive direct enumeration based on \cite{redelmeijer1981polyonimoes}. As the series grows rapidly, 
the most recent \cite{luther2011highDlatticeAnimals} has only got to 18 terms. 
The techniques in this paper have generally not been used for this or other three dimensional problems, as the number of
intermediate states is large, and if one is not sufficiently careful, can be larger than the number of objects being enumerated.

\begin{figure}
	\centering
	\includegraphics[scale=1.0]{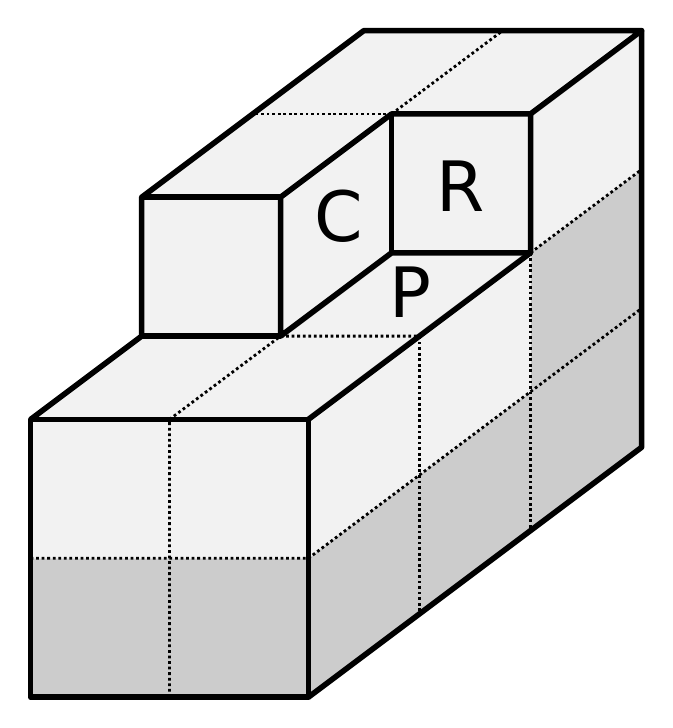}
	\caption{The boundary of a polycube finite lattice enumeration. This is the bottom of a finite lattice two wide, three deep, and at least three tall. Sites are processed starting with the
		bottom back left cube. Subsequent cubes are processed left to right. When the right boundary is reached (quickly in this case), the row starts again one step closer to the front.
		When the whole layer is done, the next layer starts again at the back left. The cubes shown are processed. The six light grey cubes are on the boundary; the nine dark grey cubes
		are behind the boundary and only matter in so far as they affect connectivity of the boundary cubes and which edge faces have been touched at least once. The next to be processed
		is the one in the top middle right indentation, adjacent to three light grey cubes labelled P (prior plane), R (prior row) and C (prior column).}
	\label{fig:PolycubeBoundary}
\end{figure}

Here is presented a transfer matrix method as a demonstration of the difficulty but not uselessness of TM on three dimensional problems. 
It is not clean, and will be a dreadful algorithm
asymptotically, but may be reasonable for intermediate values. Enumeration will be done on finite lattices, adding one site at a time, with a two dimensional
cross section as shown in figure \ref{fig:PolycubeBoundary}. 
Initial evidence indicates it may be possible\footnote{Work in progress.} to get another couple of terms using this algorithm with comparable computational
resources to \cite{luther2011highDlatticeAnimals}.

The two dimensional equivalent, polyonimoes, has been well studied and enumerated using the TM method in \cite{conway1995polyonimo}, significantly 
improved with pruning in \cite{jensen2001polyonimo}, tweaked by Knuth to add a symmetry (public talk), and parallelized in \cite{jensen2003polyonimo}.

\section{Techniques common to both algorithms}

The algorithms are very similar in the sense that the number of states is the primary driver of efficiency and memory use. The efficiency of these techniques comes
from the conflation of different states when it can be recognised that they produce the same result. The more that states can be conflated, the more efficient the algorithm.

One can think of starting with a simple algorithm that builds up the objects being enumerated one element at a time, with the state being the entire history (as in figure \ref{fig:FullDirAnTree}).
There is no gain from dynamic programming so far as these states are never reached twice. Then one works out a way of abstracting some reduced information about a state called a {\it signature}
which is sufficient to identify the state sufficiently precisely to compute downstream computations (as in figure \ref{fig:PartialDirAnTree}). 
The expectation is that many states map to one signature. One then redefines the signatures
to be the states and one has a more efficient algorithm.

So a signature is just a state, although with the connotation of being a summary of many other original states.

\subsection{Signature design}

Design of a signature is usually the most important part of the algorithm. The more states that map into a signature, the more efficient the algorithm.
The specifics depend very much on the problem, but there are a couple of approaches that can be tried.

One approach is to think about different ways of describing the state. This can be exemplified by 1324 PAPs. In \cite{conway20151324PAPs} the same basic state is used 
as in the earlier work \cite{johansson20141324PAPs}, but a different way of describing it is used. The information is basically what numbers in the permutation
are left and what restrictions there are on their use. In \cite{johansson20141324PAPs} these restrictions are, for each available number, an index pointing back at which prior
numbers are unlocked after all numbers after this one are dealt with. In \cite{conway20151324PAPs} it is noticed that this can be represented by a 
set of nested brackets. Looking at things this way, it becomes clear that consecutive brackets can be simplified, reducing the number of
states. It also makes some state factorizations more obvious.

Another approach is to think about storing different information in the state. The following sections give several examples of this.

\subsubsection{Signature design on a finite lattice}

For the frequent case of enumerating objects on a finite lattice, where each operation consists of considering one extra site to be processed, the
state generally consists of the boundary of the processed sites. This is generally one dimension smaller than the lattice being processed. The state
then consists of the elements on the boundary, plus any needed connection information. The boundary consists of the set of points in the processed
set that have at least one neighbour in the unprocessed set. Connectivity information is necessary in order to prevent counting disconnected or
otherwise invalid objects. 

In the 3D polycube case, this connectivity information is basically which sites on the boundary are connected to which other sites.
This can be practically implemented by assigning an (arbitrary) integer to each site. Sites with the same integer are connected. It is important
to {\it canonicalize} such information to prevent the same state being referred to in two different ways by different choices of
integers. A simple canonicalization in this case is to order the sites, and define the canonical choice of integers to be the one
that would come first in lexicographic order if the integers were written out as a series in site order.

For 2D lattices, the connectivity information is simpler as there are frequently restrictions on crossings. So for the 2D version
of polyonimoes, the boundary is a line (possibly with a kink in it). It is impossible to have a series of sites A,B,C,D where A 
is connected to C, and B is connected to D, but A is not connected to B. This means that connectivity information can often be
defined as ``connected to the next appropriate thing'' or ``connected to the last appropriate thing'' \cite{conway1995polyonimo}. This does not
intrinsically improve the algorithm, but it does improve implementation as it removes the need for canonicalization, and means
that the state can be stored in a small number of bits,
which is good for speed and memory usage of the big hash map typically used.

For the 2D directed animals case, the connectivity information is trivial... everything on the line perpendicular to the preferred direction
must have come from the root, so is all connected, so no connectivity information is needed, just the state of the sites
on the boundary.

For the enumeration of objects where bonds between sites have a meaning, it is worth considering the boundary along either sites
or bonds - usually one will be significantly more efficient than the other.

\begin{figure}
	\centering
	\includegraphics[scale=2.0]{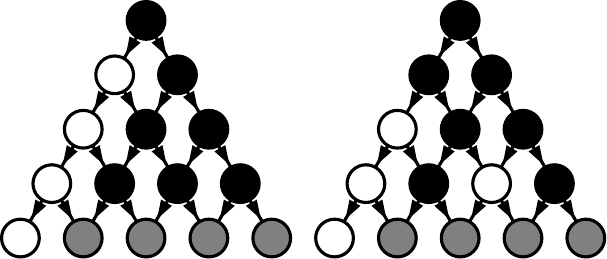}
	\caption{Two different starts of a directed animal that have different last sites occupied (the second last row), but the same sites occupiable (the last row). Black sites
		are occupied; grey sites are occupiable.}
	\label{fig:DirAnAB}
\end{figure}

Another common technique is to move the boundary forward one step to be just in front of the processed sites. That is, the sites
that are in the unprocessed set but have a neighbour in the processed set. Instead of
recording what is on the processed sites, one records what could be in the unprocessed site. Directed animals are an
excellent example. Consider the following two states when recording what is present on the boundary. State $A$ has three consecutive
sites occupied - subsequent growth can come from any of the four sites in the next row (figure \ref{fig:DirAnAB}, left). State $B$ is the same, except the middle
of those three sites is not occupied (figure \ref{fig:DirAnAB}, right). But subsequent growth for state $B$ is the same as state $A$, as the same four sites are
reachable from state A. So $f(A)=f(B)$, and $A$ and $B$ should be considered the same state. This can be done by storing which
sites in the next row could be occupied - the new state replacing $A$ and $B$ would be four occupiable sites in a row.

A similar technique can be used with 3D polycubes. The kink in the boundary will usually (depending on 
which site one is currently processing) have a single place where there are three sites adjacent to the next-to-be-added site. These three sites are
in the prior plane $P$, the prior row $R$, and the prior column $C$ (see figure \ref{fig:PolycubeBoundary}). Sites $R$ and $C$ have multiple sites 
still to be considered that they
are adjacent to, but $P$ is only adjacent to the site about to be processed. If $P$ is occupied, and connected to at least one of
$R$ or $C$, then the $P$'s presence does not affect anything, as any site added will be connected to the group containing $P$ by
dint or $R$ and/or $C$. So the number of states could be reduced slightly by removing $P$ in such a circumstance. This is a
minor optimization as it only affects one site.

Instead of storing what is connected to what in the past, one could store the information about what should be connected to what.
This is different information with a significantly different set of child states.
In \cite{clisby2012sap} this reduced the number of states slightly, but more importantly made the trimming (see later) much faster.

\subsubsection{Touching boundaries on a finite lattice}

Suppose one wants to enumerate polycubes of size 15. This can be done by enumerating all polycubes in a 15 by 15 by 15 finite
lattice. Of course smaller polycubes will appear multiple times differing only by translation. 
This can be accounted for by enumerating polycubes on all
subsets of the 15 by 15 by 15 lattice and doing appropriate subtractions\footnote{Define $U_L$ to be the total number of objects
on a lattice $L$ (or a generating function of them), and $C_L$ to be a canonicalized version only counting objects that fit into
$L$ but not any smaller lattice. Then $C_L$ can be generally computed as $U_L$ minus some multiples of $C_l$ for lattices $l$ smaller
than $L$. For the zero size lattice, typically $C_l=0$. Then by induction given all $U_l$ for $l$ a subset of $L$, one can compute $C_l$.
One can then add up all $C_l$ to get the number of objects that fit on lattice $L$ without multiple counting. In practice it is 
usually simpler than this due to constraints on touching boundaries}.

There are ways of improving this. For many symmetric lattices the result is independent of the order of the dimensions. For instance,
the number of 3D polycubes on a 3 by 5 by 7 lattice is the same as on a 7 by 3 by 5 lattice. This clearly reduces the number of lattices
that must be enumerated upon, but, more importantly, one can choose the largest dimension to be in the direction perpendicular to the
boundary. The size of the border (which is often the main driver of algorithm complexity) is determined by all but that dimension, 
so choosing it as the longest reduces the maximum size of the border. 

This means the biggest border that one will come across to enumerate polycubes of size 15 will be 5 by 6, a much more tractable
problem.

Furthermore, if one is counting polycubes on an $A$ by $j$ lattice, where $A$ is the boundary area and $j$ the depth, 
then instead of enumerating separately
for each $j$, one can just work up to the maximum $j$ and one will pass through all lower values in the process, and
whenever an object is finished in a layer $l\leq j$, one stores for the $A$ by $l$ lattice. One generally also wants to add the
requirement that, after the first layer (of size $A$) is done, the empty boundary is not accepted. This provides uniqueness in
the $j$ dimension, avoiding having to explicitly canonicalize, but, more importantly, it reduces the number of states that will
be processed, improving efficiency.

Another way to get translational uniqueness is to require each edge of the lattice
to be touched. Intuitively, this seems a bad idea, as it makes the state more complex as the state will now
contain flags indicating whether each lattice edge\footnote{Not including the two edges dealt with in the prior paragraph} has been touched.
This increases the total number of legal states, which is an upper bound on the number of states that will be reached, and a
useful heuristic for algorithm complexity.
However, in practice \cite{jensen2001polyonimo} adding the boundary requirement 
often significantly reduces the number of states that {\it will} be reached, as it will take too many sites
or bonds to get to that state and still be in a position to finish the object (see next section on trimming).

\subsection{Signature trimming}

Trimming is the process of making sure that the algorithm is clean. That is, it never produces a state $Z$ such that $f(Z)=0$.
The number of states processed, which is the main driver of efficiency, is bounded above for a clean algorithm by the number
of objects being enumerated times the maximum length of a chain from the start state to an end state (which is typically linear 
in the length of the sequence being computed). An unclean algorithm can be much worse. So trimming is of comparable importance to
state design.

An example is in two dimensional polyonimo enumeration where the TM algorithm with trimming in \cite{jensen2001polyonimo} is much more efficient 
than the prior algorithm \cite{conway1995polyonimo}, even though the states are more complex, as they involve keeping track of whether the boundaries have been reached.
Similarly it has been very effective in 2D polygon enumeration with \cite{jensen1999sap} improving on \cite{enting1980sap}, and in
many other problems.

The trimming algorithm is used on each state processed to discard useless next possible states. Sometimes this is obvious, and
not even worth mentioning. In other cases it is quite complex \cite{jensen1999sap}, which can be problematic if it is time consuming, as it is executed
for each state. In \cite{clisby2012sap} trimming is made faster by changing the state definitions from describing the connectivity in 
the processed space to describing the required connectivity in the unprocessed states. 

In the finite lattice case, trimming basically consists of determining the minimum number of sites or bonds needed to
complete the object, usually by resolving connectivity issues including to the boundaries. If this minimum number is too high
the state can be rejected. Typically ``too high'' means that the lowest term of the generating function associated
with the state plus the minimum completion number is greater than the desired length of the series. 

Sometimes it is difficult to come up with a perfect trimming algorithm, in which case an imperfect algorithm may be used which
leaves the algorithm still not clean, but better than with no trimming algorithm. An imperfect (or conservative) trimming 
algorithm generally produces a lower bound on the number of elements needed to finish
rather than the number itself. The enumeration algorithm will then still produce the correct answer,
but will waste time dealing with states that will go to zero. This is the case for the 3D polycube
algorithm presented here as an example of a difficult trim. 

\subsubsection{Polycube trimming}

Specific details of the algorithm for trimming polycubes is presented here. It is not provided as an example
of a good algorithm; to the contrary, it is a horrible algorithm. It is slow, complex (and therefore
error prone), and imperfect. Rather it is presented to show how this is difficult, and the why it is
sometimes worth changing the signature design to make trimming easier. Bad as it is, without it the
algorithm would be totally impractical. Hopefully a reader will be able to improve it!

The finite lattice method applied to 3D polycubes adds one site at a time. So there are three nested loops, the outermost
iterating over layer ($z$) then row ($y$) and then column ($x$).  The boundary is a slice through the
finite lattice, with a kink - some sites are in the layer currently being processed ($z$), some are in the prior layer
($z-1$, or empty if $z=0$).
The state consists of the occupied sites on this border, and their connectivity. This is defined as a number
for each site. Zero means unoccupied; a positive integer means occupied and connected to all other sites with the
same number. For the rest of this section, these numbers will be called {\it colors}. 
Flags are also kept for which of the four sides of the finite lattice are attached (the $z=0$ side is connected by
construction - no zero states are allowed after the last site of the $z=0$ layer is processed, and so the last $z$ layer
touched is always the one currently being processed).

The task for the trimmer is to compute what the minimum number of sites is that needs to be added to finish the polycube.
This requires adding cubes to:
\begin{itemize}
	\item connect each distinct color,
	\item connect to the unconnected sides,
	\item connect to the minimum value of $z$ that ending is allowed\footnote{By symmetry, one can always have $z$ be the biggest dimension. So polyonomoes can be required to reach 
		to a certain minimum $z$, being the maximum of the width and height}.
\end{itemize}

A fast precise answer is unknown to the author. A reasonably fast, conservative algorithm is presented here.

Ignore for the moment everything other than connecting all the colors. The only sites that need to be considered for
this are the ones one layer beyond the boundary, with an extra row above the kink to allow connectivity between above and below
the kink.
Any minimal connection using other sites could be done equivalently using just
this set of sites. Call this set of sites the {\it grid} $G$.

For any color $c$, we can define $table(c)$ to be a number for each site in $G$, being the minimum number of sites
needed to get to that cell from any cell of color $c$ on the boundary. This can be computed in polynomial time using
depth first search. Alternatively, define a consistency function $C(t)$ which takes a table $t$ and makes it consistent,
that is, reduces any value to no more than one more than any of its neighbours. Then $table(c)$ is $C$ of the table which
is infinite other than neighbours of $c$ which are $1$.

In the one color case, the color connectivity cost is trivially zero. 

In the two color ($a$ and $b$) case, the cost is the
minimum in $table(a)+table(b)-1$, alternatively the minimum cost in $table(a)$ of anything adjacent to $b$.
Indeed, the minimum cost of connecting $a$ and $b$, going via a given grid element, is $bitable(a,b)=C(table(a)+table(b))$

In the three color ($a$, $b$, and $c$) case the connectivity cost is the minimum value of $bitable(a,b)+table(c)-1$,
alternatively the minimum cost in $bitable(a,b)$ of anything adjacent to $c$.

Indeed, we can now define $tritable_c(a,b,c)=C(bitable(a,b)+table(c)-1)$. Then $tritable(a,b,c)=\min{\left\lbrace tritable_c(a,b,c),tritable_b(a,c,b),tritable_a(b,c,a)\right\rbrace}$
is now the cost of connecting $a$, $b$ and $c$ going via a particular site.

This lets us now do the four color case ($a$,$b$,$c$,$d$), in which case the connectivity cost is the minimum
value of $tritable(a,b,c)$ adjacent to $d$.

For the five or more color case, the cost of continuing this is becoming prohibitive. A conservative lower bound
is used, being just the cost to connect four of the colors\footnote{In practice this was done by computing
$tritable(a,b,c)$ and then finding for each remaining color the minimum value adjacent to it, and taking the
maximum of these values}.

Now consider the cost to connect to the edges as well. Ignoring connection costs, one could find the sites with
minimum and maximum $x$ and $y$ values, and take the distance to the edges from them. The connection cost will
be that, probably with some extra as one also usually\footnote{The exception is when an edge site is in the kink} 
need to go up a site in order to build out. If one adds in the color connection cost, frequently the sites used in that
will provide the extra height. Working out exactly when extra sites are needed is difficult. The sites
in the color connection cost will almost never\footnote{Exception: sometimes the connection cost can be done in
the row below the kink at the same cost as the row above the kink.} detour outside the minimum and maximum $x$ and $y$ values,
so the color connection cost can be safely added to this edge connection cost\footnote{With a 1 discount for the 
case of the previous footnote, when one needs to get to the bottom and the kink is up to the last or second last row.} to get a lower bound.

Similarly, to get to the maximum $z$ one has to have sites going from the current layer up to the minimum $z$ layer.
One of these may be counted in the color connection cost if it is non-zero. So this distance can be added, with a 1 discount
in the case of a non-zero color connection cost.

This ends up with a somewhat conservative trimming function which produces an algorithm that performs vastly better than
with no trimming. Asymptotically it will be a dreadful algorithm as the number of states that end up being zero will
become huge, but appears to be reasonable for currently calculable polycube sizes.
	
\subsection{Binning}

A process, sometimes described as binning, can be occasionally used to reduce memory use at the cost
of execution time. One chooses some intermediate set of states, divides them into groups, and, for each group,
one does the enumeration repeatedly, requiring that group to be passed through. Then one adds up the results. This is rarely useful
as the time penalty is usually significantly greater than the memory advantage. 

\section{Techniques particular to dynamic programming} \label{sec:TechDP}

Dynamic programming is generally inferior to a transfer matrix algorithm, but is still used extensively as it is
conceptually easier, and has some optimizations not available to the transfer matrix method.

\subsection{Factorizations} \label{sec:Factorizations}

Sometimes it is possible to factorize a state, that is say $f(S)=f(S_1)f(S_2)$, for some states $S_1$ and $S_2$.
Generally $S_1$ and $S_2$ will be significantly smaller than $S$, vastly reducing the work needed to compute them.
For instance, in the 1324 PAPs, the state $S=(S_1)S_2$ has this property, as the brackets mean that all things
to the right of the bracket must be processed before anything inside the bracket\footnote{There is another
similar but not quite as effective optimization used... if $S=S_1(S_2)S_3$ then $f(S)=\sum_i m_i f(s_i S_2)$ where
integers $m_i$ and prefix states $s_i$ are functions of $S_1$ and $S_3$ but not $S_2$.}.

As far as dynamic programming is concerned, multiplication is insignificantly different from addition. But
multiplication of two states does not fit into the transfer matrix paradigm of cumulative sums at all. 
This ability is the main advantage of the dynamic programming method over the transfer matrix method.

It is possible
to use a hybrid algorithm, where one basically uses a transfer matrix technique, but when a factorization is encountered,
the simpler state is evaluated using a dynamic programming algorithm, and then becomes a constant multiplier for
the more complex state. A constant multiplier is then fine from a transfer matrix perspective, apart from frequently having skipped 
ahead some steps in the transfer matrix. This skip makes one have to store it and merge it back in when other states
have caught up. The advance storage of states plus the extra overhead of dynamic programming storage as well somewhat undoes the
smaller memory use advantage of the transfer matrix method. It turned out to look promising but not actually
help significantly for 1324 PAPs (unpublished work). Such hybrid algorithms are almost twice as complex to 
implement as either base algorithm, increasing the likelihood of programming errors.

\subsection{Probabilistic caching}

When memory is a greater issue than speed, it is possible to reduce memory use by only caching results
probabilistically. That is, after computing $f(S)$, instead of storing $S\rightarrow f(S)$ in some table,
only store it with some probability $p$. If $f(S)$ is never needed again (as often happens), then there
is no cost in not storing it. If $f(S)$ is needed frequently, then eventually it will be stored, and there will
be no subsequent penalty. A smaller $p$ produces more memory savings, but a greater time penalty, allowing some
tuning of the algorithm to just fit in the memory available. A high quality random number generator is not
needed; one fast and simple approach is to use a simple accumulator. Each time one wants to see if something
should be stored, add $p$ to the accumulator. If the result is at least one, subtract one and store the value.

Superficially, this may sound as if it reduces the memory use to a factor of $p$, and increases the time by a
factor of $1/p$, but it is not that simple. The total number of calls made will be increased, and
so the amount stored will be $p$ times a larger number. So memory is reduced to a factor between $p$ and $1$.
Execution time increase is not as bad as may be expected because the number of times each particular state is 
referenced changes the effect of $p$ on the running time, with large and small numbers both improving matters.
It is difficult to theoretically determine this factor, as it depends on the call graph, so empirical results
are needed.

\begin{figure}
	\centering
    \includegraphics[width=\linewidth]{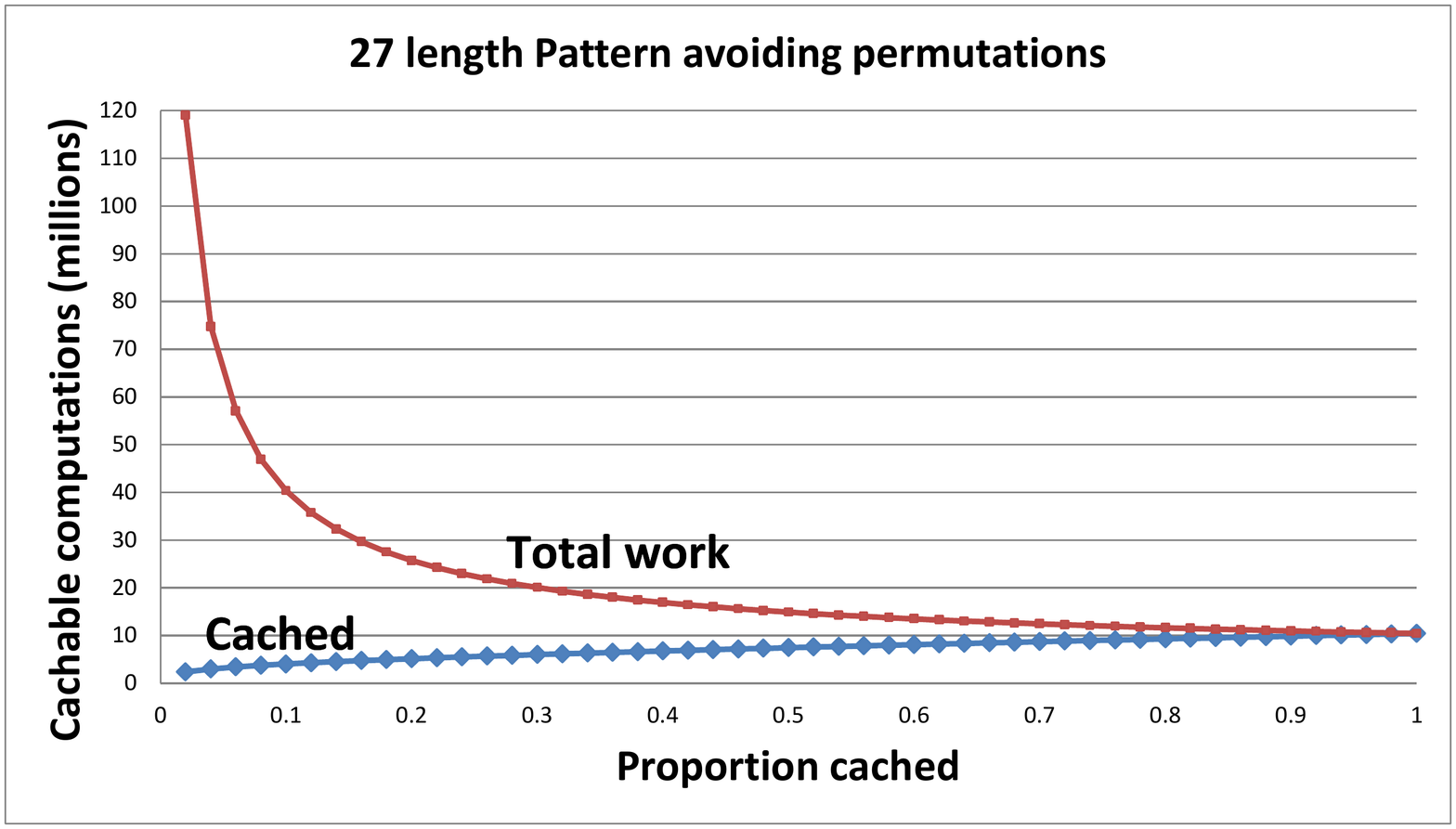}
	\caption{Probabilistic caching in 27 step PAPs as a function of the probability of caching a newly computed result. Memory use is roughly proportional to the bottom line; time is roughly proportional to the top line.}
	\label{fig:PC1324}
\end{figure}

\begin{figure}
	\centering
	\includegraphics[width=\linewidth]{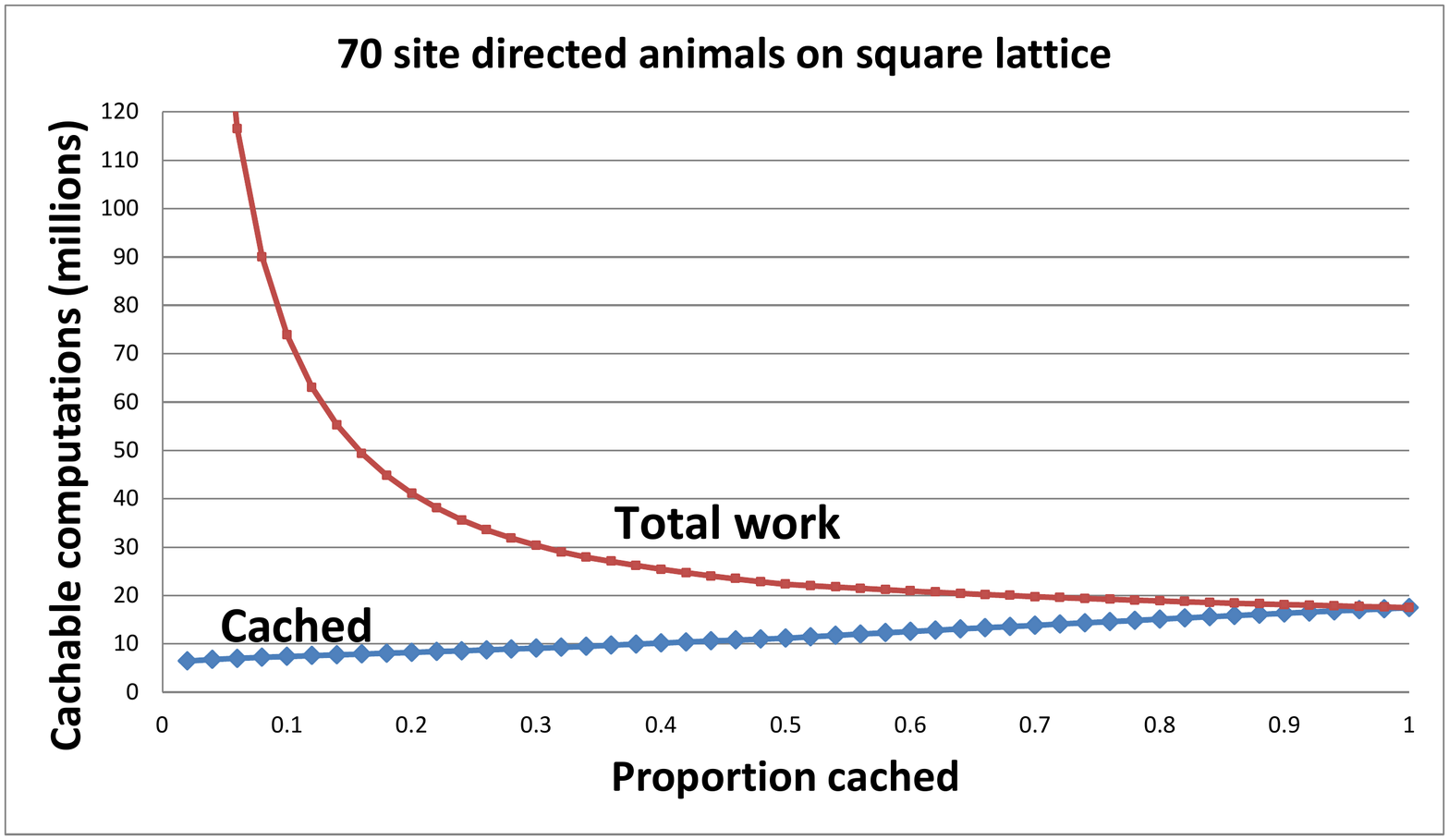}
	\caption{Probabilistic caching in 70 site directed animals on the square lattice as a function of the probability of caching a newly computed result. Memory use is roughly proportional to the bottom line; time is roughly proportional to the top line.}
	\label{fig:PCDirAn}
\end{figure}

Empirical results show that this works surprisingly well 
and was used in the memory constrained enumeration of 1324 PAPs \cite{conway20151324PAPs} to get an extra term. Figures
\ref{fig:PC1324} and \ref{fig:PCDirAn} show the empirical effects of $p$ 
on directed animals and 1324 PAPs\footnote{The
actual time and memory use for the 1324 PAPs is a little more complex due to factorizations.}. 
The patterns are surprisingly similar given that the directed animals have a maximum
of two children while the 1324 PAPs can have dozens. A $p$ of $0.3$ gives roughly a 40 percent reduction
in memory and takes roughly twice as long in the two cases shown.

This optimization cannot be used in transfer matrix algorithms as the data store is a cumulative sum, rather than
a cache.
 

\subsection{Cache inspection}

One is trying to canonicalize states such that two states that can be proven to produce the same result are actually encoded as one state.
In the absence of infinite wisdom,
a method that can be used to check that nothing obvious has been missed is to inspect the cache after the algorithm has run to see if there are multiple states
with the same value. To reduce the chance of coincidences, one can run the algorithm until there are a hundred thousand or so entries
in the cache, and look at those entries with associated values over a million. If there are multiple states producing the same large
value, an inspection of those states will hopefully inspire a realisation of some equivalence of states, producing a more
efficient definition of states.

This has been recently applied by the author to the 1324 PAPs algorithm in \cite{conway20151324PAPs} to notice
that signatures starting with a 1 or a comma allow commutation of some of the ending terms. Proving this then led to other useful realisations
leading to a more efficient algorithm (not yet published).

This is easy to do for a dynamic programming algorithm; for a TM algorithm a similar technique can be used, but is less directly 
associated with state equivalence. 
	
\section{Techniques particular to transfer matrix}

There are also some techniques primarily applicable to the transfer matrix method.

\subsection{Part of the signature is the iteration number} \label{sec:TMItNoSig}

For directed animals, the obvious state is the combination of the sites on the boundary eligible for occupancy, and
the number of sites left to include. For children there are two main implementation choices:
\begin{itemize}
	\item One can take a processing step to consider a certain eligible site on the boundary, and have
	one child if it is occupied, and another if it is not occupied. In the first case the
	number of sites left to include will have decreased, in the second the number of available sites on the boundary will
	have reduced. One can define a hierarchy function based primarily upon the number of sites left to include and
	secondarily on the number of available sites; this hierarchy function will be lower for both child states,
	so loops are impossible. The call graph for this algorithm is shown in figure \ref{fig:PartialDirAnTree}.
	\item One can take a processing step to mean the next used site is chosen. There can be many substates
	in this case (one for each available site in the boundary). The number of sites left to include is
	decreased by one each time; it makes an ideal hierarchy function; loops are impossible. 
	The call graph for this algorithm is shown in figure \ref{fig:PartialDirAnTree2}.
\end{itemize}
For dynamic programming implementations the two are comparable. The latter is slightly more complex to implement
but will be slightly faster and use a little less memory. However for a transfer matrix style implementation
they are significantly different. Both algorithms can be implemented with each iteration processing one
value of the hierarchy function. For the first algorithm this is somewhat fiddly and has the issue that
a child may have a value of the hierarchy function a few steps further on; these need to be tracked.
For the second algorithm, each child will go directly into the list to be processed at the next iteration
as the hierarchy function is ideal.

In both cases, the number of sites to be included is implicit in the current TM iteration number. Therefore
one does not have to store it as part of the signature, saving memory in the states, a large number
of which will be stored. This is shown in table \ref{tab:TMDirAn} where the third argument of
each state is determined by the row number in the table, and therefore does not need to be stored with each state.

This implicit storage of part of the state is obvious for the finite lattice methods where each iteration
corresponds to moving the boundary out by one site. The shape of the boundary is identical for
each element of a particular iteration, and so does not have to be stored with each state.

\subsection{Signature Invariants} \label{sec:siginvariant}

Sometimes it is possible to divide the states into groups that have the property that there is no state that is a child
of a member of two different groups (or bins). For instance, in the 3D polycubes case, define a group by the occupied/unoccupied
status of each site other than the one about to be covered. Adding (or not) a new site may affect connectivity, but it
won't change other sites' occupancy. These are signature invariants

This division means that instead of processing each signature sequentially, putting all the results in a giant hash map,
one can do each group separately, using a small hash map (with better cache locality), and then, when a group is finished,
extract the states and associated multiples, store them in some more compact format, and reuse the hash table for the
next group. This enables the use of a sparsely filled hash map, improving speed, without the massive memory hit of having
a huge sparse hash map. One does have to be very careful with the overhead of small groups and cache clearing.

Processing this way, it is often possible to make the construction of the next set of groups implicit. In the 3D
polycube case, it makes sense to assign each group a number whose binary representation contains
the occupied status of each site, with the most recent site added as the most significant bit. Then the output 
from the processing of
each group will be in one of two groups, identified by shifting the current groups' number down one bit and adding
a new 0 or 1 as the most significant bit. These can be serialized into two piles based on the most significant bit,
and then when all groups are processed, the two piles are concatenated. If the groups started off in order, they
will now again be in order, this time for the next set of groups. Since the initial null state is by definition in order,
the groups will always be in order.

This can be used to efficiently use disk as storage instead of memory. Current low latency SSDs are
too slow to be used as swap space for a giant hash map (2016 unpublished tests), but using them as storage for the 
serialized group outputs requires many fewer random accesses and is reasonable in some situations.

Of probably greatest importance, this can be very useful on multiprocessor systems, as each group can be assigned to a node,
and the node can process that group knowing it will not have to share the hash map with any other node. Examples of such use
include 2D self avoiding walks \cite{conway1996saw,jensen2013saw}, polygons \cite{jensen2003parallelsap} and polyonimoes \cite{jensen2003polyonimo}.

This sort of construction is often possible for geometric entities on finite lattices; it is less clear for
other cases like 1324 PAPs which has proven difficult to parallelize.

\subsection{Finite lattice techniques}

Finite lattice enumerations generally involve enumerating objects on a set of finite lattices individually,
and reconstructing the total number of objects from the results on finite lattices. 

Finite lattice enumerations are ideal for the TM algorithm as each iteration is well defined (add one more
portion, typically site, of the lattice). The techniques described here can be used for DP algorithms as
well, but there is rarely any reason to use DP rather than TM for such algorithms.

\subsubsection{Finite lattice symmetries}

As mentioned before, when enumerating all the objects on a symmetric lattice of a certain size, the result is often
identical to enumerating it on other sizes determined by the lattice symmetry. For instance,
enumerating 3D polycubes on a 3*8*2 lattice is the same as a 2*3*8 lattice (and four others).
This reduces the number of finite lattices that need to be computed. 

This is not always possible - if one of the sides of the lattice is special (e.g enumeration
of paths that are attached to one side of the lattice) or if the dimensions are different
(e.g. enumeration of polygons on the square lattice by both horizontal and vertical bonds).

Generally it is better to make the boundary go across the shorter dimension(s), as this
makes the states simpler, and probably less numerous, although with good trimming this can
be less important than one might expect.

\subsubsection{State reflection}

When the iteration is such that the kink does not preclude symmetry (e.g. when finishing a row
or a plane), frequently the boundary condition is equivalent to its mirror image. At this point
a consolidation can be done, declaring one of these arbitrarily to be the canonical one, and
merging the two. 

This can in principle halve the number of states to process, although this is somewhat misleading
as the number of states will continue to grow up to almost the number it would have been anyway
during the subsequent iterations where the kink prevents this consolidation. It also breaks
many of the invariants (section \ref{sec:siginvariant}).

\subsubsection{Reconstruction from fragments}

Some of the objects being enumerated can be split up into fragments and then reconstructed
by enumerating those fragments. This can reduce the number of items being enumerated.
The fragments being enumerated are often called {\it irreducible}.

For instance, a bridge on a lattice $L$ is a path on $L$ that has one end at one side of $L$,
and the other at the opposite side of $L$. The bridge is irreducible if there is no plane slicing
through L parallel to the attached sides of $L$ that intersects only one bond. For many problems
one can reconstruct all bridges from just the reducible bridges by chaining them with single bonds
connecting the irreducible bridges. That is, if $b(x)$ is the bridge generating function, and
$b_i(x)$ is the irreducible bridge generating function, then 
$$b(x)=b_i(x)+xb_i(x)^2+x^2b_i(x)^3+\dots=\frac{b_i(x)}{1-xb_i(x)}$$

Enumerating irreducible objects can potentially be faster than enumerating all objects as
there are fewer of them, and therefore there will probably be fewer states 
reached. Also this can be used to reduce the size of the signature as irreducible
objects are more compact and fit onto narrower lattices. This was critical to \cite{conway1993saw} but 
later advances in understanding of pruning made this less useful.

\section{Associated values}

The associated values for dynamic programming usually are integers as that fits in well with
the paradigm. The state is then the boundary condition and number of elements yet to be included. 

With the transfer
matrix, one could do the same thing, although it is often more efficient to have the state
be the boundary condition, and the associated value is then an array of integers; the polynomial
coefficients of a generating function multiple. 
This is more efficient as the state only needs to be stored once,
and trimming calculations only need to be done once for each state. This is difficult to
do with a dynamic programming algorithm as one usually does not know in advance how big
the generating function will need to be, which is essential information to pass to the
child state evaluator.

The number of non-zero entries in a generating function is usually quite small. This is
because most of the states are complex, as there are only small numbers of simple states,
and complex states typically take a lot of elements to produce, leaving
few spare elements for the series. To take advantage of this, avoid storing zero elements. 
In the majority of cases, the non-zero elements are consecutive, so that means
each generating function can be represented by a start index, a length, and an array
of coefficients of that length.

\begin{figure}
	\centering
	\includegraphics[width=\linewidth]{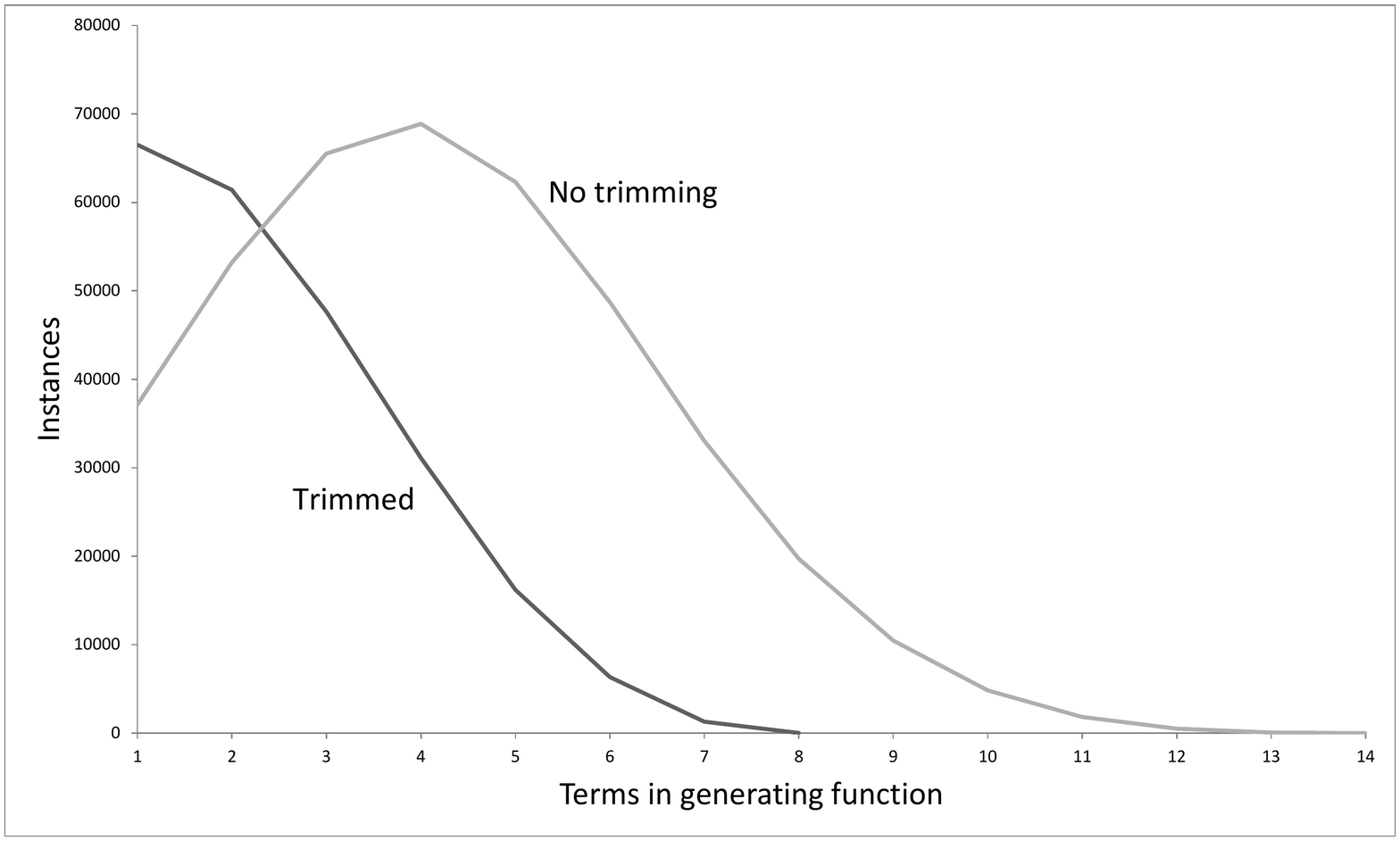}
	\caption{Number of generating functions of each length encountered midway through the
      enumerating of up to 16 element polycubes on and spanning the 3 by 5 by 5 lattice. Computations
      were done with and without using the (imperfect) trimming.}
	\label{fig:GFLength}
\end{figure}

Of course if most generating functions are only length 1, then that is a lot of overhead to
store one integer, but the total overhead is generally less storing generating functions, especially
if the alternative of storing the number of elements left with the state just makes it have
similar extra overhead instead.
See figure \ref{fig:GFLength} for an example distribution of lengths, and the effect of
even imperfect trimming on reducing not just the total number of states, but the size of
the associated generating functions as well.
 
Sometimes only even (or odd) terms can be non-zero, in which case only those terms should be stored.
This is particularly common with enumeration by bonds.

\subsection{Multi variable series and moments}

It is straight forward to have more than one variable in the generating function. Suppose instead of
counting polycubes just by number of cubes, one also wants to count by surface area (3D perimeter). 
Let $A_{i,j}$ be the number of polycubes with $i$ cubes and perimeter $j$. Instead of having 
a 1D array as the generating function associated with each state, one would have a 2D generating
function, and update it in the same way. Trimming becomes a little more complex to describe, as
does the highest term of the enumeration.

Sometimes however, that can take up too much memory, and having a moment series is useful. 
The $m$-th moment of surface area series, counting by volume, would then be 
$$M_m(x) = \sum_{i,j} A_{i,j}x^i j^m$$
Note that $M_0(x)$ is just the normal enumeration by volume (number of cubes).

In a finite lattice computation, when a state is being processed for a site, the new
generating function needs to be modified by adding $c$ cubes and surface area $s$.
With a 2D generating function, this is straightforward, one just shifts it $c$ in one
dimension and $s$ in another dimension. For $M_0(x)$ it is also straightforward,
just shift $c$. For other moments it is not quite as obvious, but it turns out
to be straightforward. Let $M^*_m$ be the new desired moment, and $M_m$ be the 
existing moment. Then 
$$M^*_m(x) = \sum_{i,j} A_{i,j}x^{i+c} (j+s)^m = x^c \sum_{i,j} A_{i,j}x^i \sum_{k=0}^m \binom{m}{k} j^{m-k}s^k =  x^c \sum_{k=0}^m \binom{m}{k}s^k M_{m-k}$$
which means if one has existing moments from $0$ to $m$, one can easily compute the new 
moments from $0$ to $m$.

Given moments, one can compute, say, the mean surface area for size $i$ polycubes
by dividing the  $i$-th coefficient of $M_1$ by the $i$-th coefficient of $M_0$.
Different variables allow different properties of the object to be studied.

This is described in detail in \cite{conway1995polyonimo} for the 2D polyonimo case.
By a similar process one can compute percolation series, where the generating function
is or the form $\sum A_{i,j} p^i(1-p)^j$ storing just a single generating function for each state.
	
\section{Implementation issues}

Generally, an efficient algorithm is more powerful than an efficient implementation.
An efficient algorithm can be many orders of magnitude better than a competitor, whereas
an efficient implementation is typically only a couple of orders of magnitude better than
a simpler, less efficient implementation. With that said, the implementation does matter.

\subsection{Chinese remainder theorem}

The numbers being computed often are larger than the native size of an integer. Many languages
have libraries for big integers (data structures representing large integers). However these
tend to be slow and memory consuming.

One can use the Chinese remainder
theorem to deal with this problem. One does the enumeration with small integers by doing everything
modulo some number (typically a prime). Do this for several different, coprime, moduli, and this 
is enough information to be able to simply regenerate the number
modulo the product of the original moduli. This is generally faster, and more memory efficient
than using a big integer library. It also has the advantage of offering some redundancy; a
concern in long running computations can be an error somewhere in the computer (e.g. from
gamma ray strikes). Such an error in one of the moduli will usually produce a very large,
noticeable, inconsistency in the results.

This is probably the most commonly used trick in implementation.

Knuth\footnote{Public talk}, implementing TM algorithms, faced with the inefficiency of repeating the computations
multiple times, separated out the state part from the value part. He made two programs, the
first of which would carry around all the state information, and would produce as output
a long stream of instructions for another program to execute, to actually do the additions and trimming of generating functions.
This is more complex and introduces significant disk IO, but has two advantages:
\begin{itemize}
	\item The first program only needs to be run once, while the second program can be
	run for each modulus, rather than recomputing the states for each modulus. This is particularly
	useful if trimming is computationally expensive and therefore a majority of the time.
	\item The memory consumption can be slightly lower (in principle up to a factor of 2) than doing
	both at the same time.
\end{itemize}

\subsection{Parallelization}

Multi-processor machines are the norm, and to use a big computer effectively, a parallel
algorithm is usually needed. This is usually quite difficult, as distributed hash maps
take a large amount of inter-processor communication which can be the bottleneck,
even on shared memory systems due to the expense of cache synchronization.

Most parallelized success stories with transfer matrices have used signature invariants
(see section \ref{sec:siginvariant}). If speed rather than memory is the constraint,
then the problem can be partitioned cleanly by different moduli and lattice sizes.

Writing parallelized code is very difficult and error prone. Writing efficient
communication code is even more difficult and requires a deep understanding of
the architecture used. One can use distributed libraries, although they generally
don't do exactly what one wants and have an associated overhead. There are no
clear general answers.

TM algorithms are usually easier to parallelize than DP algorithms, as the distributed hash
map lookup of DP is particularly unfriendly to parallel architectures.

\subsection{Hash map}

Both TM and DP require a large amount of storage for a map from state to value. For
DP this is the cache; it must be a map as the operations are {\em look up value for a state} and
{\em store calculated value for a state}. For TM the operations are {\em iterate over state,value pairs} and
{\em store value for a state, adding to existing if available}. These are most obviously performed
by a map, although the map-reduce framework (such as Apache Hadoop) is also possible. Also if
invariants are used, a small hash and long list may be used. But most frequently a map (usually a hash map) is used.

This hash map is usually very large. This makes some of its properties affect the performance
of the program significantly. The load factor of the hash map is important - higher values 
mean a slower program using less memory. 

Having a good hash function of the state is very important as
the state tends to have many similar values, and a poor hash function can lead to a very
large number of collisions.

Most languages come with support for hash maps built into a standard library. However this
is generally optimized for convenience and flexibility rather than performance, and a less
general purpose one can give order of magnitude improvements particularly to memory but also time. 

For instance, in Java, the standard hash map is from one object to another object. It is
often possible to encode the state as a bit string and therefore as an integer. Similarly
the values are often integers\footnote{Integer is used as a generic name for a natively 
handled data type. It may be called long or something other than integer.}. When the
generic map is used from integers to integers, then each integer has to be wrapped in
an object, causing a large overhead. Similarly, many standard library hash maps
handle collisions by storing at each entry a list of values that mapped to that entry.
This causes another level of wrapping for each entry. There exist much more memory efficient
(and fast) hash maps, such as the {\em GNU Trove} library which includes multiple versions for each
combination of primitives as key or value. 

It is frequently useful to have multiple maps instead of one big one. Some function of
the key determines which map to use, and some other function of the key is used in
that map. Reasons for doing this include:
\begin{itemize}
	\item Many hash map libraries have limited size hash maps. Many use 32 bit integers
	in the indexing (and hash values). On current computers this can be a serious limit
	to the size of the map.
	\item If one doesn't set the size of a map in advance, it will typically automatically
	increase in size when needed. This operation is usually implemented by creating a new
	data store of larger size and copying values across. During this process both
	the new and old datastore are in memory, which can exceed available memory. By
	splitting up the maps, this situation causes less transient memory increase.
	\item In a multithreaded, shared memory program, it can reduce contention for write locks.
	\item Some of the key can be implicitly stored in the map selection stage. In \cite{conway20151324PAPs}
	the keys for 1324 enumeration were encoded as 128 bit integers. The first 64 bits were 
	used to determine which hash map to use (in a reversible manner). The remaining 64 bits
	were used as the key in the map. This reduced the memory used for the keys in the maps.
\end{itemize}

\begin{figure}
	\centering
	\includegraphics[width=\linewidth]{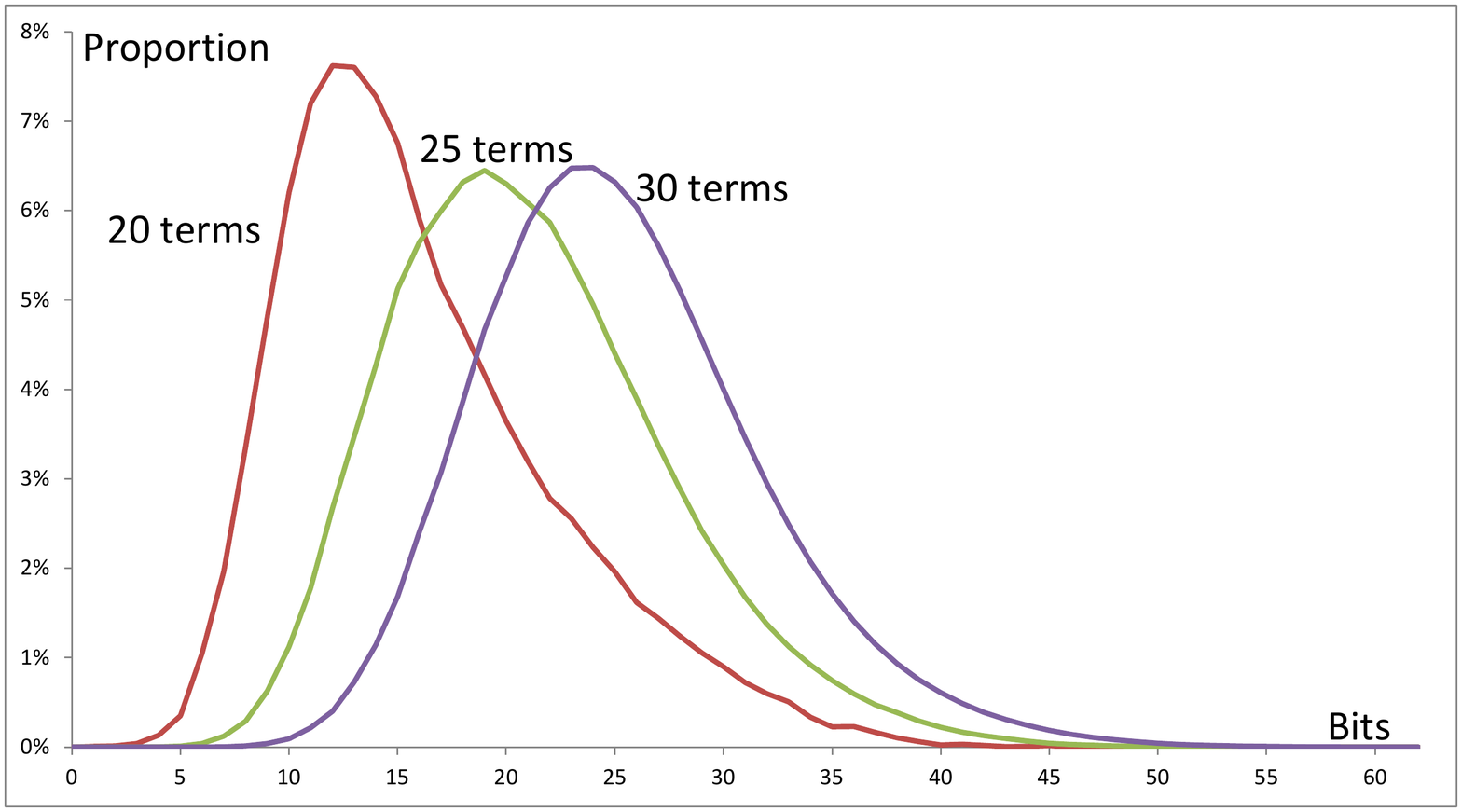}
	\caption{Bits used in cache for 1324 PAPs, length 20, 25 and 30, as a proportion of the total cache entries. The full answer for length 30 is 78 bits}
	\label{fig:Bits1324}
\end{figure}

Another way of reducing memory use at the cost of time is by using multiple hash
maps with different size values depending upon the actual size of the value. Suppose the
values are 64 bit integers. Many of the actual values may fit in a 32 bit integer.
This is because a large portion of the states produced will have very restricted paths
to endpoints resulting in modest actual values. This means one uses two hash maps,
one mapping to 64 bit values and one to 32 bit values. Reads have to check both maps;
writes just go to the best fitting map. This is easy for DP; for TM it may require
changing which map a value is stored in when another value is added to it. Of course
one could get an even greater effect by using 32 bit moduli in the Chinese remainder
theorem, but the multiple maps method can achieve much the same effect with a smaller
performance hit. Profiling of the bit length in the 1324 PAPs algorithm found it looked
useful to have 24 bit, 40 bit, and 64 bit results (see figure \ref{fig:Bits1324})

\section{Conclusion}

The art of efficient enumeration algorithms is, like most skills, enhanced by knowledge of
a host of techniques, only a few of which will be appropriate for any specific problem. This
paper has described many such techniques with consistent terminology and a discussion of their
relative merits and applicability. This will hopefully help the reader to realise
when given tricks can be used. The exercise of writing this paper certainly helped me
realise how I could have improved various algorithms I came up with in the past.

In particular, the transfer matrix technique is well known and used for two dimensional lattice enumeration problems where
it seems natural, but it is rarely used outside of that domain. However, totally unrelated
problems with an ideal hierarchy function can often use the
transfer matrix technique instead of the dynamic programming that has
been typically used. This can provide a reduction in memory usage and
significant advantages for parallelization. This paper has emphasised the similarities
of the two techniques and the issues in converting from dynamic programming to transfer matrix.

Each enumeration problem will have its own special requirements and peculiarities, but many 
tricks can be reused.

Special thanks to Tony Guttmann for introducing me to enumeration, being a fantastic supervisor, and many comments on this manuscript.
	
\bibliographystyle{hunsrt}
\bibliography{combinatorics} 
	
\end{document}